
\input amstex
\documentstyle{amsppt}
\magnification=1200
\parindent 20 pt
\NoBlackBoxes

\define \im{\imath}
\define \cG{\Cal G}
\define \cH {\Cal H}
\define \Mod{\operatorname{Mod}}
\define \MU {\Mod_{U(\cG),f}}
\define \MUh {\Mod_{U(\cG)[[h]],f}}
\define \End{\operatorname{End}}
\define \Pent{\operatorname{Pent}}

\define \Alt{\operatorname{Alt}}

\define \mdh {\quad\text{mod}\quad h}

\define \cO {{\Cal O}(G)}

\define \cC{\Cal C}
\define \cN{\Cal N}
\define \cR{\Cal R}
\define \la {\langle}
\define \ra {\rangle}
\define \lm {\lambda}
\define \r{\rightarrow}
\define \lr{\longrightarrow}

\define \De{\Delta}
\define \de{\delta}
\define\as {A^\star}
\define \ot{\otimes}
\baselineskip 20pt
\vskip .20in
\topmatter
\centerline{\bf Cohomological construction of quantized}
\centerline{\bf universal enveloping algebras}
\centerline{Joseph Donin and Steven Shnider\footnote{Supported by a
grant from the Israel Science Foundation administered by the
Israel Academy of Sciences and Humanities}}
\centerline{Department of Mathematics, Bar Ilan University}
\abstract
{Given an associative algebra $A$, and the  category, $\cC$,
 of its finite dimensional modules,
 additional structures on the algebra $A$ induce corresponding ones
 on the category $\cC$.
Thus, the structure of a rigid quasi-tensor (braided monoidal)
 category on $Rep_A$ is induced by
an algebra homomorphism $A\to A\otimes A$ (comultiplication),
coassociative up to conjugation by
$\Phi\in A^{\otimes 3}$ (associativity constraint) and
 cocommutative up to conjugation by  $\cR\in A^{\otimes 2}$
(commutativity constraint), together with an antiautomorphism
(antipode), $S$, of $A$
 satisfying the certain compatibility conditions. A morphism of quasi-tensor
 structures  is given by an element $F\in A^{\otimes 2}$ with suitable
induced actions on $\Phi$, $\cR$ and $S$. Drinfeld defined
such a structure on $A=U(\cG)[[h]]$ for any semisimple Lie algebra
$\cG$ with the usual comultiplication and antipode but
nontrivial $\cR$ and $\Phi$ and proved that the
corresponding quasi-tensor category is isomomorphic
to the category of representations of the Drinfeld-Jimbo (DJ)
quantum universal enveloping algebra (QUE), $U_h(\cG)$.

   In the paper we  give a direct cohomological construction of
the  $F$ which reduces $\Phi$ to the trivial associativity constraint,
without any assumption on the prior existence of a strictly
coassociative QUE. Thus we get a new  approach to  the DJ
quantization.  We prove that $F$ can be chosen to satisfy some
additional invariance conditions under (anti)automorphisms
of $U(\cG)[[h]]$, in particular, $F$ gives
an isomorphism of  rigid quasi-tensor categories. Moreover, we prove
that for  pure imaginary values of the deformation parameter, the
elements $F$, $R$ and $\Phi$ can be chosen to be formal
unitary operators on the second and third tensor powers of
the regular representation of the Lie group associated
to $\cG$ with $\Phi$ depending only on even powers of the
deformation parameter. In addition,
we consider some extra properties of these elements and give
their interpretation in terms of
additional structures on the relevant categories.}
\endabstract
\endtopmatter
\subheading{\S0 Introduction}

The original interest in quantum groups and quantized universal
enveloping algebras (QUE) came from the existence of a universal
$R$-matrix satisfying the Yang-Baxter equation. From the
 point of view of representation
theory  the Yang-Baxter equation is a consequence of the
more fundamental  hexagon identities, which are the basic identities
for the braiding in a braided monoidal (quasi-tensor) category. However from
the categorical perspective  there is
something unnatural in the standard axioms for quantum
groups since there is  a second basic operator, the associativity
constraint, which is restricted to be the identity. A more
natural, although technically more involved, set of axioms
was introduced by Drinfeld in \cite{Dr1, Dr2} for the
algebraic structure which he called  ``quasi-triangular
quasi-Hopf algebra.'' The axioms for such structures are
precisely those required in order that the category of finite
dimensional representations be a  rigid braided monoidal category,
without  restricting  the associativity constraint to be the identity.

We  review some basic definitions.
A quasi-bialgebra is an associative algebra, $A$, with  comultiplication
$\Delta:A \rightarrow A\otimes A$, which
is not necessarily coassociative. However there is given  an invertible
element of the triple tensor product, $\Phi\in A^{\otimes 3}$,
which expresses the relation
between the two iterated comultiplications  by the formula
 $$\Phi(\De \otimes id)\De(a)=(id\otimes \De)\De(a)\Phi
\quad\text{for}\quad a\in A.\tag1$$
An additional condition on  $\Phi$ is  the identity in $A^{\otimes 4},$
$$(id^{\ot 2}\ot \De)(\Phi)\cdot(\De\ot id^{\ot 2})(\Phi)=
(1\ot\Phi)\cdot(id\ot\De\ot id)(\Phi)\cdot(\Phi\ot 1).
 \tag2$$
This equation comes from requiring the pentagon identity in the
category of finite dimensional representations of $A$ as
explained in the next section.
In a  quasi-triangular  quasi-bialgebra we are given a so-called  $R$-matrix,
 $\cR$, which is an invertible element of the double tensor
product, $A^{\otimes 2}$,  expressing  the
relation between the comultiplication
and the opposite comultiplication, $\De^{\text {op}}=\sigma\circ\De$  for
$\sigma(a\ot b)=b\ot a,$ by the formula
$$ \cR \De(a) =\De^{\text {op}}(a) \cR,
\quad\text{for}\quad a\in A.\tag3$$
The compatibility conditions between $\cR$ and $\Phi$  are given by the
identities in $A^{\otimes 3}$,
$$
\align
(\De\ot id)\cR&=\Phi^{312}\cR^{13}(\Phi^{132})^{-1}\cR^{23}\Phi \tag{4a}\\
(id\ot \De)\cR&=(\Phi^{231})^{-1}\cR^{13}\Phi^{213}\cR^{12}\Phi^{-1}.\tag{4b}
\endalign
$$
These two equations correspond to two commutative hexagons diagrams
which, according to the MacLane coherence theorem,
 together with the pentagon diagram corresponding to (3)
generate all the relations involving  the associativity constraint and the
commutativity constraint in a quasi-tensor category.
A quasi-bialgebra with an antipode (see \S 1)  is called a quasi-Hopf algebra.

In the context of deformation theory, Drinfeld proved two important results
about the existence and uniqueness of such structures.
Let $\cG$ be a Lie algebra over a field $K$ of characteristic 0,
$K[[h]]$ be the ring of formal power series with
coefficients in $K$, and  $U(\cG)[[h]]=U(\cG)\otimes_K K[[h]],$
the completed tensor product.
The first result says that for any symmetric $\cG$-invariant element $t\in
\cG^{\otimes 2}$ there exists a  $\cG$-invariant
element  $\Phi_h\in U(\cG)^{\otimes 3}[[h]]$
satisfying the pentagon identity and such that together
with $\cR_h= e^{\pi \im ht}$ it satisfies the hexagon identities.
These elements define a deformation (quantization) of
the universal enveloping algebra as a quasi-triangular quasi-Hopf
algebra, $(U(\cG)[[h]],\cR_h,\Phi_h)$.  The basic algebraic operations,
multiplication, comultiplication and antipode, are undeformed, being
defined on $U(\cG)[[h]]$  by the $K[[h]]$
linear extension of the standard operations on $U(\cG)$.
The second result says that,
 modulo an equivalence relation corresponding  to an
equivalence of braided monoidal categories, this deformation is unique.

A non-trivial Hopf algebra deformation of
$U(\cG)$, the Drinfeld-Jimbo quantum universal enveloping algebra (DJ
QUE),  was constructed by explicit
formulae about 1985 \cite {Dr3}\cite {J}.  The   existence of such
a deformation together with the
uniqueness theorem just mentioned proves that
 Drinfeld's quasi-Hopf deformation  has an equivalent Hopf
presentation. The equivalence is defined by an element
$F_h\in(U(\cG)^{\otimes 2}[[h]]$ which transforms $\Phi_h$ to $1$
and conjugates the comultiplication, see equations (17) and (26)
below.  From this point of view the existence $F$  follows from
the prior  knowledge of the existence of a Hopf deformation.

In our approach the existence of $F$  is proved
directly. This gives an new,  less ad hoc, construction of the
DJ QUE and suggests an approach to the construction
of other  examples.
We use the methods of the classical deformation theory of algebras as
developed by Gerstenhaber,  with Cartier coalgebra cohomology replacing
Hochschild algebra cohomology. In a kind of secondary obstruction
theory we shall also use Chevalley-Eilenberg  Lie algebra
cohomology.

The paper is organized as follows:
In \S 1 we give some of the basic definitions from category theory
and explain the axioms for the antipode in a quasi-Hopf algebra.
In \S2 we review Drinfeld's cohomological proof of the existence
of $\Phi$ and show how to include some additional symmetries,
the significance of which have been explained in \S1. Then in
\S3 we prove our main result, which says that,
given an infinitesimal $f\in\wedge^2\cG$, which induces the structure
of a Lie bialgebra on $\cG$, in order for there to exist a
 Hopf (as opposed to quasi-Hopf) quantization of
$U(\cG)$, it is enough that a certain subcomplex of  the Chevalley-Eilenberg
complex of $\cG^*_f$ (the induced structure on the dual of the
Lie algebra $\cG$) have zero cohomology in dimension 3.
In \S4 we apply this result to the DJ infinitesimal for
a simple Lie algebra and give a purely cohomological proof
of the existence and uniqueness up to equivalence of the
DJ quantum group. Our proof is in much the
same spirit as the cohomological proof of the existence of
quantizations of certain Poisson brackets, \cite{DS,Li,DL}.
Having  dealt in detail with the construction of the associativity constraint,
$\Phi$, in \S2, in the appendix we sketch
 Drinfeld's cohomological proof of the existence of
a pair  $(\Phi_h,\cR_h)$ satisfying the pentagon and hexagon
identities.  We prove that for pure imaginary deformation
parameter, $\bar h =-h,$ the elements $F,R,$ and $\Phi$ can be
chosen to be formal unitary operators on the second and third
tensor powers of the regular representation of the Lie group
associated to $\cG$. Moreover, we consider some extra properties
of these elements and give their interpretation in terms of additional
structures on the category of representations.

\subheading{\S1 Preliminary categorical remarks}

Recall that a {\bf monoidal category} is a triple $(\cC,\ot,\phi)$
where ${\Cal C}$ is a category equipped with a functor $\ot:\cC\times\cC\to
\cC$, called the tensor product, and a functorial isomorphism
$\phi:(M\ot N)\ot P\to M\ot (N\ot P)$
called associativity constraint.  The latter satisfies the pentagon
identity, that is, the diagram
\NoBlackBoxes
$$
\CD
((M\ot N)\ot P)\ot U @>\phi>>
(M\ot N)\ot(P\ot U) @>\phi>>
 M\ot (N\ot (P\ot U))  \\
@V{\phi\ot id}VV @.  @A{id\ot\phi}AA\\
(M\ot (N\ot P))\ot U @.\overset\phi\to
{-\!\!\!\!-\!\!\!\!-\!\!\!\!-\!\!\!\!-\!\!\!\!-\!\!\!\!-\!\!\!\!
-\!\!\!\!-\!\!\!\!-\!\!\!\!-\!\!\!\!-\!\!\!\!-\!\!\!\!-\!\!\!\!
-\!\!\!\!-\!\!\!\!-\!\!\!\!-\!\!\!\!-\!\!\!\!-\!\!\!\!-\!\!\!\!
-\!\!\!\!-\!\!\!\!-\!\!\!\!-\!\!\!\!-\!\!\!\!-\!\!\!\!-\!\!\!\!}@>>>   M\ot
((N\ot P)\ot U)
\endCD\tag5
$$
is commutative. In addition we assume the existence of an object
{\bf 1} which is a two sided identity for $\ot$ and such that the
composition
$$M\ot N\lr (M\ot {\bold 1})\ot N\overset \phi \to \lr M\ot
({\bold 1}\ot N)\lr M\ot N\tag6$$
is the identity.

Let $A$ be a Hopf algebra over a field $K$,
of characteristic zero, and let $\Mod_A$ be the
category of $A$ modules which are finite dimensional vector spaces
over $K$. First we define the structure of a
monoidal category on $\Mod_A$ using the comultiplication.
 Given two representations, $\rho_M:A\r \End(M)$ and
$\rho_N:A\r \End(N)$ on
$M$ and $N$ respectively, the tensor product (over $K$) of vector spaces,
$M\otimes N$ is naturally a representation of the tensor product of algebras
$$\rho_M\otimes\rho_N: A\otimes A\r \End(M)\otimes \End(N)
\cong \End(M\otimes N).$$
 Composition with the  comultiplication
$$\rho_{M\otimes N}=(\rho_M\otimes \rho_N)\circ\De\tag7$$
defines a representation of $A$ on $M\otimes N$.
Since the comultiplication is coassociative we
can take as the associativity constraint in  $\cC$ the identity morphism.
The element ${\bold 1}$ is given by the field $K$ with $A$ module
structure coming from  the augmentation of $A$,
$\epsilon:A\r K.$
For a quasi-Hopf algebra we use the same definition of the
tensor product, but since the usual coassociativity condition
is replaced by (1), the associativity constraint is given by the
action of $\Phi$ arising from the natural $A^{\ot 3}$ module structure on
 the triple tensor product of $A$ modules. Equation (2) implies
the commutativity of the pentagon (5). In order to guarantee
condition (6), $\Phi$ must satisfy
$$(id\ot \epsilon \ot id)\Phi=1.$$

We say that a monoidal category $\cC$ has a
{\bf rigid } monoidal structure if for each object $M$
there is a ``left dual" object $M^*$ and a ``right dual" object
$^*\!M$ together with morphisms
$$ \text{(a) }\quad {\bold 1}\r M\ot M^*,\quad\text{(b)}
\quad  M^*\ot M\r {\bold 1}, \quad\text{(c)}\quad  {\bold 1}\r  ^*\!\!M\ot
M,\quad\text{(d)}\quad M\ot ^*\!\!M\r {\bold 1}.\tag8$$
For the left dual object we require that the following
two  compositions give  identity morphisms,
$$M\r {\bold 1}\ot M\r (M\ot M^* )\ot M \overset \phi \to
\r M\ot (M^* \ot M)\r M\ot{\bold 1}\r M,\tag{9a}$$
$$M^*\r M^*\ot{\bold 1}\r M^*\ot(M\ot M^* ) \overset \phi^{-1} \to
\r (M^*\ot M) \ot M^*\r {\bold 1}\ot M^*\r M^*,\tag{9b}$$
and similar diagrams for the right dual.

The antipode of a (coassociative) Hopf algebra, $A$, is an operator, $S:A\r A$
satisfying,
$$m (S\ot id) \De(a)=\epsilon(a)=m(id\ot\Delta) S(a).$$
A rigid monoidal structure on the category $\Mod_A$ is given by
defining  the left dual as the vector space dual with (left module)
action given by $a\cdot \lambda=\lambda\circ S(a)$. The $A$ module morphisms,
$ {\bold 1}\r M\ot M^*$ and $M^*\ot M\r {\bold 1}$, in (8a,b)
are defined by $k\mapsto k\cdot id$, where
$M\ot M^*$ is identified with $\End(M)$, and $\lambda\ot x\mapsto \lambda(x)$
respectively.
The vector space structure on the right dual is the same
as the left dual but the module structure is given by
$a\cdot \lambda=\lambda\circ S^{-1}(a).$

In the case of modules  over  a quasi-Hopf algebra, conditions
(9a,b) give two equations relating the antipode and $\Phi$.
We introduce two new elements $\alpha, \beta\in A$
coming from the definition of the morphisms in (8).
 $${\bold 1}\r M\ot M^*\quad\text{is given by}\quad
k\r k\cdot(\rho_M(\beta)\ot 1)\circ id$$
 and
$$ M^*\ot M\r {\bold 1}\quad\text{is given by}\quad
\lambda\ot x\r \lambda(\rho_M(\alpha) x).$$
Then the  axioms for the antipode in
a quasi-Hopf algebra, (equations 1.17-1.19 of \cite{Dr1}) and
equation (10) below
guarantee that $\Mod_A$ is a rigid monoidal category.
We use the Sweedler notation $\De(a)=\sum a_{(1)}\ot a_{(2)}$
and sometimes delete the summation sign for simplicity in notation.
$$\align
\sum a_{(1)}\beta S(a_{(2)})&=\epsilon(a)\beta\quad
\text{corresponding to diagram 8(a)}\tag{10a}\\
\sum S(a_{(1)})\alpha a_{(2)}&=\epsilon(a)\alpha\quad
\text{corresponding to diagram 8(b)}\tag{10b}\\
\sum \Phi_1\beta S(\Phi_2)\alpha \Phi_3&=1\quad
\text{corresponding to diagram 9(a)}\tag{10c}\\
\sum S((\Phi^{-1})_1)\alpha(\Phi^{-1})_2\beta S((\Phi^{-1})_3)
&=1\quad\text{corresponding to diagram 9(b)}.\tag{10d}
\endalign
$$

We will be interested in two supplementary conditions on $\Phi$,
$$\align\Phi^{321}\Phi&=1\quad\text{and}\tag{11}\\
\Phi^{321}&=\Phi^S,\tag{12}
\endalign$$
where the superscript $S$ indicates applying the antipode to
all three tensor components.

Condition (11) is a particular case (for $\cR=1$) of the
equation
$$ \cR^{21}(\Delta\otimes id)\cR =\Phi^{321}[\cR^{23}
(id\otimes\Delta)\cR]\Phi\tag{13}$$
which appears in the definition of what Drinfeld  calls a coboundary
structure.
Condition (12) is equivalent to the compatibility of the associativity
constraint with the rigid structure as expressed by the
commutativity of the diagram

$$\CD
[M\ot (N\ot P)]^*@>(\phi_{M,N,P})^*>> [(M\ot N)\ot P]^* \\
@V \cong VV @V\cong VV\\
(P^*\ot N^*)\ot M^* @>\phi_{P^*,N^*, M^*}>> P^*\ot(N^*\ot M^*)
\endCD\tag{14}
$$
where $\cong$ are compositions of the natural equivalence
$(M\ot N)^*\r N^*\ot M^*.$ We shall return to the interpretation of equations
(11) and  (13)
after discussing the concept of equivalence of monoidal categories.

\define\wt{\widetilde}
Let $(\cC,\ot,\phi)$  and
$(\wt{\cC},\wt{\ot},\wt{\phi})$ be two monoidal categories,
then a {\bf monoidal functor} from $\cC$ to $\wt{\cC}$ is given by a pair
$(\chi, \eta)$ where $\chi:\cC\to \wt{\cC}$ is a functor
and $\eta:\chi(M\ot N)\to \chi(M)\wt{\ot}\chi(N)$ is a functorial
isomorphism such that the diagram

\define\al{\chi}

$$
\CD
\al ((M\ot N)\ot P) @>\eta>>
\al (M\ot\ N)\wt{\ot}\al (P)
@>{\eta\wt{\ot}id}>>(\al (M)\wt{\ot}\al (N))\wt{\ot}\al (P)\\
@V{\al(\phi)}VV @. @V{\wt{\phi}}VV\\
\al (M\ot (N\ot P)) @>\eta>> \al (M)\wt{\ot}\al (N\ot P)
@>{id\wt{\ot}\eta}>>\al (M)\wt{\ot}(\al (N)\wt{\ot}\al (P))
\endCD\tag{15}
$$
is commutative and $\al({\bold 1}_{\cC})\cong {\bold 1}_{\cC'}.$

For a functor between  rigid monoidal categories we will
require the additional conditions, $\chi(M^*)\cong \chi(M)^*$
and $\chi(^*\!M)\cong ^*\!\chi(M)$, as well as commutativity
 of the diagram

$$\CD
\chi((M\ot N)^*)@<\cong<<(\chi(M\ot N))^*
@<(\eta_{M,N})^*<< (\chi(M)\tilde\ot\chi(N))^* \\
\quad@V\cong VV @.@ V\cong VV\\
\chi( N^*\ot M^*) @>\eta_{N^*, M^*}>>\chi(N^*)\tilde\ot\chi(M^*)
@>\cong>>\chi(N)^*\ot \chi(M)^*
\endCD\tag{16}
$$
and of a similar diagram for the left dual.

Given a quasi-Hopf algebra $(A, m, \De, S, \Phi)$,
let  $F$ be an element of  $A^{\otimes 2}$ and
define a transformation of the comutiplication on $A$ by
 $$\tilde\De=F\De F^{-1}.\tag{17}$$
Relative to this new comultiplication there is
 a new tensor product structure on $\Mod_A$  defined by
$$\rho_{M\tilde\otimes N}=
F_{M,N}\circ\rho_{M\otimes N} \circ F_{M,N}^{-1},\quad
\text{where}\quad F_{M,N}=(\rho_M\otimes \rho_N)(F).\tag{18}$$
Then diagram (15)  defines the associativity constraint,
$\tilde\phi$, for $\tilde\ot$.
It is given by representing the element
$$
\tilde\Phi=(1\otimes F)((id\otimes \De)F)\Phi
((\De\otimes id)F^{-1})(F^{-1}\otimes 1).\tag{19}$$

In other words, $\chi=id$ and
$$\eta_{M,N}=F_{M,N}:M\ot N\lr M\tilde\ot N,\tag20$$
define an equivalence of
monoidal categories $(\Mod_A,\ot,\phi)$
and $(\Mod_A,\tilde\ot,\tilde\phi)$.
If in addition
$$(F^{21})^S F=1\tag21$$
diagram (16) and its left dual counterpart are  commutative and
we have an equivalence of rigid monoidal categories.

In terms of quasi-Hopf algebras, transforming $\De$ by (17),
$\Phi$ by (19), leaving the antipode unchanged, and transforming the elements
$\alpha$ and $\beta$, respectively, to
$$\tilde\alpha=\sum S(F^{-1}_{1i})\alpha F^{-1}_{2i},
\quad\text{and}\quad\tilde\beta= \sum F_{1i}\beta S(F_{2i})\tag{22}$$
defines a transformation {\bf twisting}, which, in turn,
 determines an equivalence relation.
If two quasi-Hopf algebras,  $A$ and $A'$, are equivalent under
twisting then the
rigid monoidal categories $\Mod_A$ and $\Mod_{A'}$ are equivalent.

Consider $A=U(\cG)$, the universal enveloping
algebra of a Lie algebra $\cG$, and $\MU$, the category of
finite dimensional representations.
 We define  a  deformation of the monoidal category $\MU$
by first extending the coefficients from the field $K$ to
the formal power series algebra $K[[h]]$.
For any finite dimensional module $M$ we define the free $K[[h]]$
module of finite rank,
$M[[h]]=M\otimes_K K[[h]]$. Let $\MUh$ be the category
consisting of $U(\cG)[[h]]$ modules which are free
$K[[h]]$ modules of finite rank.
Define  comultiplication on $U(\cG)[[h]]$ as the
$K[[h]]$ linear extension of comultiplication on $U(\cG)$.
Together with the tensor product of $K[[h]]$ modules this
 defines a monoidal structure on $\MUh$ with associativity
constraint given by the identity. Relative to this monoidal
structure the imbedding defined above defines
 an imbedding of monoidal categories.

Drinfeld's quasi-Hopf deformation is related to a
nontrivial  monoidal structure on $\MUh$.
The tensor product of modules is the standard one given by (7) with the
usual comultiplication  but the associativity operator will be nonstandard:
$$\phi_{M,N,P,h}=(\rho_M\otimes \rho_N\otimes \rho_P)\Phi_h,\tag{23}$$
where $\Phi_h$ is a $\cG$ invariant element in
$U(\cG)^{\otimes 3}[[h]]$ satisfying the pentagon
identity (2).

Let $(\cC,\ot,\phi)= (\MUh,\otimes,\phi_h)$
and $(\tilde{\cC},\tilde\ot,\tilde\phi_h)= (\MUh,\tilde\otimes,\phi_h)$
where $M\tilde\ot N=N\ot M$ and the associativity constraint is
the natural one, $\tilde\phi_{M,N,P,h}=\phi_{P,N,M,h}^{-1}.$
Then equation (11) for $\Phi_h$ is equivalent to the
commutative diagram (15) where the natural transformation
$\eta$ is given by transposition.
Equation (13) is a   generalization of equation (11)  in which
the  functor $\chi$ is the identity and
the natural transformation between $M\ot N$ and $M\tilde{\ot}N=N\ot M$
is given by composing transposition with the $\cR$-matrix.

Using the standard
antipode we can also define a rigid monoidal structure.
In this case, the elements $\alpha$ and $\beta$ will be chosen to be
 invariant, allowing us to reorder the product, and will satisfy
$$\alpha\beta=(\sum \Phi_1 S(\Phi_2)\Phi_3)^{-1},\tag{24}$$
for example $\alpha=1,\quad\quad \beta=(\sum \Phi_1 S(\Phi_2)\Phi_3)^{-1}.$

For any Lie algebra $\cG$ defined over  a field
of characteristic zero, the  existence of such
a $\Phi_h$  satisfying (11) has been proven in \cite{Dr1,Dr2}.
In the next section we outline the proof and show that essentially the
same arguments prove that $\Phi_h$  can be chosen so that it
satisfies (12) as well. (In the case when $\cG$ is simple the
$\Phi_h$ is unique up to change of parameter and ``twisting"
as defined in (20) and (22) below. \cite{SS})

The invariance condition on $\Phi_h$ guarantees that
$\phi_h$ defines a natural transformation in $\MUh$. The pentagon
equation (2) for $\Phi_h$ guarantees the commutativity of (5).
The triple  $(\MUh,\otimes,\phi_h)$ defines one type
of deformation of the  monoidal category $(\MU,\ot, id)$.
This is not yet the deformation which gives the quantum group
since  the tensor product has not been deformed. However the
required deformation is given by an equivalent rigid monoidal
category.

As an aside, not required in any of the subsequent proofs, but
useful in understanding the situation,
we relate the deformation of the category to
deformation of the algebra $\cO$ of representative functions
generated by the matrix entries of $\rho_M$ as $M$
runs over the objects of $\MU$.
The multiplication on
$\cO$ has a simple relation to the tensor product of modules.
For any pair of elements
$m\in M$ and $\lm\in M^*$ define
$$\la u,f_{m, \lm}\ra=\la \rho_M(u)m,
 \lm \ra.$$
Then the fact that multiplication in $\cO$ is dual to comultiplication
in $U(\cG)$ implies immediately that
$$\align
\la u,f_{m,\lm}f_{n,\mu}\ra&=
\la \De(u), f_{m,\lm}\ot f_{n,\mu}\ra=\sum \la u_{(1)}, f_{m,\lm}\ra
\la u_{(2)}, f_{n,\mu}\ra\\
&=\sum\la\rho_M(u_{(1)})m, \lm\ra \la \rho_N(u_{(2)})n, \mu\ra=
\la (\rho_M\ot \rho_N)(\De(u))(m\ot n), \lm\ot\mu\ra\\
&=\la \rho_{M\ot N}(u) (m\ot n), \lm\ot\mu\ra
=\la u, f_{m\otimes n, \lm\otimes \mu}\ra.\tag{25}
\endalign$$
One  approach to understanding the quantized function algebra
is in  terms of a deformation of the tensor product  on the
monoidal category $\MU.$
Since any module $M_h$ in $\MUh$ has the form $M\otimes_K K[[h]]$,
the space of matrix coefficients  in
$\MUh$  can be identified with $\cO[[h]]=\cO\ot K[[h]]$
and a tensor product on $\MUh$ which reduces modulo $h$ to
the standard tensor product on $\MU$ defines a formal deformation
of $\cO.$ We return to this topic briefly at the end of \S 3.

In \S 3 we describe a purely cohomological proof
of the existence of an element $F$ which
transforms  $\Phi_h$ to $1$ and conjugates the comultiplication as
in (17).  We prove that, under
some very natural assumptions on the algebra $\cG$, there exists an
$F_h\in U(\cG)^{\ot 2}[[h]]$  satisfying (21) and
$$
(1\otimes F_h)((id\otimes \De)F_h)\Phi_h
((\De\otimes id)F_h^{-1})(F_h^{-1}\otimes 1)=1.\tag26$$

This $F_h$ will define an equivalent rigid monoidal category
to $(\MUh, \ot,\phi_h)$ with
associativity constraint equal to the identity. At  the level
of quasi-Hopf algebras, we have $\Phi=1$ so $\tilde\alpha\tilde\beta=1.$
Substituting $\tilde S$ defined by
$$\tilde S(a)=\tilde\beta S(a) \tilde\alpha\tag27$$
in equations (10a) and (10b) we see that
$\tilde S$ is an antipode satisfying the usual axioms
for a Hopf algebra.

\subheading{\S2 Construction of a nontrivial  associativity constraint}

In this section we prove a modified version of Drinfeld's theorem,
 \cite{Dr1}, on the existence of a nontrivial associativity constraint
compatible with the undeformed tensor product. The modification
involves the inclusion of additional symmetries.

\proclaim{Theorem 2.1}
Let $\cG$ be  a Lie algebra over a field $K$ of characteristic zero,
 $U(\cG)$ its universal enveloping algebra,
$K[[h]]$ be the algebra of formal power series with
coefficients in $K$, and $U(\cG)[[h]]$ the Hopf algebra over
$K[[h]]$ given by extending all operations $K[[h]]$ linearly.
Let $\varphi\in \wedge^3\cG$ be  $\cG$ invariant. Then there
exists a  $\cG$ invariant formal power series
$\Phi=\Phi_h\in U(\cG)^{\ot 3}[[h]]$ solving the pentagon identity (2)
and satisfying the following conditions:

(a) $\Phi=1 + h\varphi\quad\text{mod}\quad h^2.$

(b) $\Phi^{321}\Phi=1.$

(c) Let $\theta$ be any   automorphism of $\cG$
leaving $\varphi$ invariant. If we extend $\theta$ in the natural way to $
U(\cG)^{\ot 3}[[h]]$, then $\Phi$ can be chosen satisfying,
in addition to the conditions above,
$$\Phi^\theta=\Phi\tag{c}$$

(d) If we  extend the antipode $S$ to $ U(\cG)^{\ot 3}[[h]]$ in
the natural way, then $\Phi$ can
be chosen satisfying, in addition to (a),(b), and (c),
$$\Phi^S\Phi=1\tag{d}$$

Given a commuting set of $\theta$, we can find $\Phi$ invariant under
the entire set.
\endproclaim

\noindent{\bf Proof.}
Rewrite the pentagon identity in the form
$$\Pent(\Phi):=(1\ot\Phi)\cdot(id\ot\Delta\ot id)(\Phi)\cdot(\Phi\ot 1)
\cdot[(\Delta\ot id^{\ot 2})(\Phi)]^{-1}\cdot
[(id^{\otimes 2}\otimes\Delta)(\Phi)]^{-1}=1.\tag28$$
As usual in the theory of formal deformations we try to solve (28)
 recursively in powers of $h$.
Assuming we have defined $\Phi^{(n-1)}$, an $(n-1)^{st}$
 order polynomial in $h$
satisfying (28) to order $h^n$, we
look for an element $\varphi_n\in U(\cG)^{\ot 3}$ such
that $\Phi^{(n)}=\Phi^{(n-1)}+ \varphi_n h^n$ solves (28) to order
$h^{n+1}$.  The ``obstruction'' to such an extension is a cocycle in
the coalgebra cohomology  $U(\cG)$,
whose definition we shall review quickly.

The  complex has degree $n$ component
$\cC^n=U(\cG)^{\otimes n}$, and the
coboundary is induced by imbedding $\cC^n$
in the Hochschild $n$ cochains on $\cO$  with values in $K$
via the pairing of $U(\cG)$ and $\cO$.
$$\la u_1\otimes \cdots \otimes u_n, f_1\otimes \cdots\otimes f_n\ra=
\la u_1, f_1\ra\cdots \la u_n, f_n\ra .\tag29$$
The pullback of the Hochschild coboundary to  $\cC^n$ is
$$\delta(u_1\otimes\cdots \otimes u_n)=1\otimes (u_1\otimes\cdots
\otimes u_n) - \De(u_1)\otimes\cdots \otimes u_n +\cdots +
(-1)^{n+1}(u_1\otimes\cdots\otimes u_n)\otimes 1,\tag{30}$$
This coboundary operator, first introduced by Cartier \cite{Ca},
defines a cohomology controlling deformations of a coassociative coalgebra
in just the same way that Hochschild cohomology controls
the deformations of associative algebras.
We are interested in the complex of $\cG$ invariants, which we denote
by $\tilde{\cC}$.
 The cohomology is well known, see \cite{Dr1}  or \cite{SS}:
$$H(\cC)=\wedge \cG \quad\text{and}\quad
H(\tilde{\cC})=(\wedge\cG)^{\cG}.$$
The obstructions to extending $\Phi^{(n-1)}$ are in
$H^4(\tilde{\cC})=(\wedge^4\cG)^{\cG}$. Since the relevant cohomology
group is in general non-zero, the
obstruction cocycle for extending a truncated solution
of (28) to one higher order may not  cobound.
Following Drinfeld, we show that condition (b) involves
restricting our attention to a  subcomplex with zero
cohomology in dimension 4. Define an involution, $\tau,$ of $\tilde{\cC}$
commuting with $\delta$ by
$$\tau(u_1\otimes\cdots \otimes u_n)=
(-1)^{\frac{n(n+1)}2}(u_n\otimes\cdots\otimes u_1).\tag31$$
Let $\cC_{\tau,\pm}$ be the
subcomplex consisting of the $\pm 1$
eigenspaces  $\tau$. The complex splits as a direct sum of the
two subcomplexes corresponding to the two eigenspaces.
The cohomology of each subcomplex is the corresponding eigenspace of
the restriction of $\tau$ to $\wedge\cG$.
Since the action of $\tau$ on $\wedge^n\cG $
 is multiplication by $(-1)^n$,
$$H^n(\cC_{\tau,+})=0\quad\text{for n odd, and}\quad
H^n(\cC_{\tau, -})=0\quad\text{for n even.}$$
This carries over to the subcomplex of $\cG$ invariants as well.

Suppose that we have satisfied condition (b)
up to order $h^n$, write
$$(\Phi^{(n-1)})^{321}\Phi^{(n-1)}=
(\Phi^{(n-1)})^\tau\Phi^{(n-1)}=1 + \eta_n h^n \mdh^{n+1}.\tag32$$
Assuming that $\Phi^{(n-1)}$ was $\cG$ invariant,
 $\eta_n$ is also. The commutativity of
$\Phi^{(n-1)}$ and $(\Phi^{(n-1)})^{\tau}$ implies
that $\eta^{\tau}_n=\eta_n$, so if we define
$$\Phi^{(n-1)'}=\Phi^{(n-1)}+\frac12 \eta_n h^n,\tag33$$
then $\Phi^{(n-1)'}$ satisfies (b) to order $h^{n+1}$
and is $\cG$ invariant.

Substituting in the pentagon identity we have
$$\Pent(\Phi^{(n-1)'})=1+ \xi_n h^n\tag34$$
for $\xi_n\in \tilde{\cC}^4.$
The fact that $\xi_n$ is a cocycle is a rather tedious calculation
given in detail in \cite{Dr1}.
Apply $\tau$ to both sides and use the fact that (b) holds
to order $h^{n+1}$, we find that
$$1+\xi_n^{\tau}h^n=(\Pent(\Phi^{(n-1)'}))^\tau=
(\Pent^*(\Phi^{(n-1)'})^\tau)=
\Pent^*(\Phi^{(n-1)'})=(\Pent(\Phi^{(n-1)'}))^{-1}=
1-\xi_n h^n\quad\text{mod}\quad h^{n+1},$$
where
$$\Pent^*(\Phi):=
[(id^{\otimes 2}\otimes\Delta)(\Phi)]^{-1}\cdot[(\Delta\ot id^{\ot
2})(\Phi)]^{-1}\cdot(\Phi\ot 1)\cdot(id\ot\Delta\ot id)(\Phi)\cdot
(1\ot\Phi).$$
Thus $\xi_n$ is a cocycle in $\tilde{\cC}^4_{\tau,-}$ but
  $H^4(\tilde{\cC}_{\tau,-})=0$. Therefore $\xi_n$ cobounds an
invariant cochain $\varphi_n\in\tilde{\cC}^3_{\tau,-}.$
Define
$$\Phi^{(n)}=\Phi^{(n-1)'}+ \varphi_n h^n.$$
Since $\varphi_n^{\tau}=-\varphi_n$ we have not disturbed (b) mod $h^{n+1}$
and we have solved the pentagon identity up to
order $h^{n+1}$.

Next we note that given $\theta$ as described in (c)
there is an obvious extension, also denoted $\theta$,
 to the complex $\tilde{\cC}$ which commutes with $\de$ and
$\tau$. Since $\theta$ is an involution, when
we consider the restriction to $\tilde{\cC}_{\tau, -}$
we have a decomposition into the direct sum of two subcomplexes,
the $+1$ eigenspace and the $-1$ eigenspace and, as before, the cohomology
also splits as a direct sum. By assumption $\varphi$ is $\theta$
invariant,  so we can consider solving our
deformation problem in the subcomplex given by the $-1$ eigenspace of
$\tau$ and the $+1$ eigenspace of $\theta$. If we assume that
the $\Phi^{(n-1)}$  defined above was $\theta$ invariant, then
since $\theta$ is an algebra automorphism, the element $\eta_n$ and
the obstruction cocycle, $\xi_n$, are
easily seen to be $\theta$ invariant and the extension to order
$h^{n+1}$  can be chosen to be $\theta$ invariant together with the
other required conditions.

Finding an extension satisfying condition (d)
involves a second preliminary correction.
Assume that $\Phi^{(n-1)'}$ has been defined as in (33) so that
condition (b)  is satisfied to order $h^{n+1}$
and that condition (d) has been satisfied to order
$h^n$. Define $\chi_n$ by
$$ (\Phi^{(n-1)'})^S\Phi^{(n-1)'}=1+\chi_n h^n\quad\text{mod}
\quad h^{n+1}.\tag35$$
Then
$$\align
1+\chi^{321}_nh^n&=(\Phi^{(n-1)'321})^S\Phi^{(n-1)'321}\\
&=((\Phi^{(n-1)'})^{-1})^S(\Phi^{(n-1)'})^{-1}\\
&=(\Phi^{(n-1)'}\Phi^{(n-1)'S})^{-1}\\
&=1-\chi_nh^n \quad\text{mod}\quad h^{n+1}.
\endalign
$$
Thus $(\chi_n)^{321}=-\chi_n$ and, as before, $\chi_n^S=\chi_n.$
Therefore if we set $\Phi^{(n-1)''}=\Phi^{(n-1)'}
+\frac12\chi_n h^n $ we satisfy do not disturb condition (b) and
satisfy condition (d) to one higher order.
Thus the preliminary corrections for conditions
(b) and (d) can be done simultaneously.
If both (b) and (d) are true to order $h^{n+1}$, then
$$ 1+(\xi^{\omega}_n +\xi_n)h^n=(1+\xi_nh^n)^{\omega}(1+ \xi_n
h^n)=\Pent^*(\Phi^{(n-1)'\omega}=\Pent(\Phi^{(n-1)'})=1\quad\text{mod}
\quad h^{n+1}
$$
where $\omega$ is either $S$ or $\tau.$
Therefore the obstruction cocycle $\xi_n$ is in the subcomplex given by the
intersection of the $-1$ eigenspaces of $\tau$ and $S$.
As before $\xi_n$ cobounds a cochain $\varphi_n$ of the same
type. The element $\Phi^{(n)}=\Phi^{(n-1)'} +\varphi_n h^n$
satisfies all the required conditions  to order $h^{(n+1)}$.\qed

Note that if we have a commuting set of automorphisms,
$\theta_i$,  then all the decompositions into eigenspaces
of $\tau$, $\theta_i$ and $S$ are  compatible and we can solve
in the subcomplex consisting of $- 1$ eigenvectors of $\tau$, $+1$
eigenvectors of $\theta_i$ and $+1$ eigenvectors of $S$.

In particular, for $\cG$ a simple Lie algebra we can let $\theta$ be the Cartan
 involution corresponding to a choice of Cartan subalgebra and  system of
positive roots and, this
proves the existence of a $\cG$ invariant $\Phi$ satisfying (11) and  (12).

Several remarks are in order before closing this section.
First of all, in the next section we want to consider another $\Phi$
which is given by replacing the deformation parameter $h$ with $h^2.$
will be necessary in order to construct the  equivalence $F$.
Henceforth, when referring to $\Phi$ we intend this form.
Second, we note that the invariance of $\Phi$ implies the invariance of
$\gamma=\sum \Phi_1 S(\Phi_2)\Phi_3\in U(\cG)[[h]]$. If
we chose $\alpha$ and $\beta$ invariant and $\alpha\beta=\gamma^{-1}$
and let $S$ be the $K[[h]]$ linear extension of
the standard antipode for $U(\cG)$, then equations (10a-d) are
satisfied, making $U(\cG)[[h]]$ a quasi-Hopf algebra.
For example, we can chose $\alpha=1$ and $\beta=\gamma^{-1}$.

Finally, the construction given above does not involve the braiding.  Drinfeld
\cite{Dr2}
has given  a cohomological proof of
the existence of a pair $(\cR, \Phi)$ satisfying (3), (4ab) and (11)
 thus defining a quasitriangular quasi-Hopf structure on
$U(\cG)[[h]]$. This proof, which is much more delicate, is sketched in  the
appendix where we  also discuss the modifications necessary to deal with
additional symmetries.

\subheading{\S 3 Construction of a quantized universal enveloping algebra}

In this section we study sufficient conditions for the
existence of a solution to equation (26).
These conditions are satisfied, in particular, in the case of the
 Drinfeld-Jimbo infinitesimal for $\cG$ a simple Lie algebra.
 Twisting the  $\Phi_h$ constructed above  by
$F=F_h\in U(\cG)^{\otimes 2}[[h]]$ and transforming the antipode to $\tilde S$
as given by (27),
defines a  quantized universal enveloping algebra.
Note that this presentation is not the standard one since the
multiplication is undeformed.

Our proof uses a combination  of  coalgebra cohomology  and
Chevalley-Eilenberg cohomology for the Lie algebra structure on $\cG^*$
associated to the leading nonconstant term of $F$.
Recall that  we are assuming that $\Phi$  has the form
$$\Phi=1 +\varphi h^2 +\sum_{m>1}\varphi_{2m}h^{2m},
\quad\text{with}\quad\varphi\in
(\wedge^3\cG)^{\cG}\quad\text{and}\quad \varphi_{2m}\in (U(\cG)^{\ot
3})^{\cG},\tag36a$$
and we set $$F=1 + fh + \sum_{n>1} f_n h^n, \quad \text{with}\quad
f,f_n\in U(\cG)^{\ot 2}.\tag36b$$
The element $f$ will be called the infinitesimal of $F$.

It will be convenient to reformulate equation (26) in the form
$$B(F_h,\Phi_h)=(1\otimes F_h)((id\otimes \De)F_h)\Phi_h
-(F_h\otimes 1)(\De\otimes id)F_h=0.\tag37$$

Consider the $h$ and $h^2$ terms in (37):
$$ B(1+ fh +f_2 h^2, 1+ \varphi h^2)= (\de f) h +
(f^{23}(f^{12}+f^{13})-f^{12}(f^{13}+f^{23}) +\de f_2
+\varphi)h^2\quad\text{mod}\quad h^3,\tag38$$
 where $\de$ is the Cartier coboundary, (30).

Two remarks:
\roster
\item
The vanishing of the first term implies that $f$
must be a $\de$ cocycle.
\item
If $\Phi$ were expanded in powers of $h$, instead of powers of $h^2$,
and $\varphi$ was the coefficient of $h$ then (37) would require
$\de f+\varphi=0$ which implies $\varphi=0$ since $\wedge\cG$ is transversal
to the subspace of coboundaries.
\endroster

If the multiplicative group,
 $1+ hU(\cG)[[h]]$ acts on the left on $U(\cG)^{\ot 2} [[h]]$ by
$$u\bullet F= (u\ot u)F \De(u^{-1}),\tag39$$
then
$$B(u\bullet F,\Phi)=(u\ot u\ot u)B(F,\Phi)(\De\ot id)(\De (u^{-1})).
\tag40$$
Therefore this action carries
 solutions of (37) into solutions and we call two
solutions equivalent if they are in the same orbit.
 If $F'=u\bullet F$ for $u=1 + u_1 h$ mod $h^2,$
$u_1\in U(\cG),$
then the infinitesimals are related by $f'=f+\de u_1.$
Since the $\de $ cohomology in dimension 2 is $\wedge^2\cG$,
in any equivalence class there is a representative with infinitesimal
 $$f\in\wedge^2\cG,\tag41$$
so making this assumption  in
(36b) involves no loss of generality.

Now consider the requirement on the pair $f,f_2$ implied by (38).
$$
f^{23}(f^{12}+f^{13})-f^{12}(f^{13}+f^{23})\tag42$$
must belong to the cohomology class of $-\varphi$.
The fact that (42)  is a cocycle follows from the following lemma.

\proclaim{Lemma 3.1}
Let $F^{(n-1)}$ be  an $(n-1)^{st}$ order polynomial in $h$
of the form (36b) with infinitesimal (41) and
suppose that $F^{(n-1)}$ satisfies (37) to order $h^n$. Define
the obstruction cochain, $\xi_n$, by
$$B(F^{(n-1)},\Phi)= \xi_n h^n\quad \text{mod}\quad h^{n+1}.\tag43$$
\roster
\item The cochain $\xi_n$  is a cocycle.
\item
If $\de f_n=\xi_n$ then $F^{(n)}=
F^{(n-1)}+ f_n h^n$ defines an extension satisfying (37) to order $h^{n+1}.$
\endroster
\endproclaim

\par\noindent{\bf Proof.} Multiplying (43) on the right by
$[\De \ot id (F^{(n-1)})(F^{(n-1)}\ot 1)]^{-1}$
and recalling the form of $B$ given in (37),
$$
\align
\tilde\Phi&=
(1\otimes F^{(n-1)})((id\otimes \De)F^{(n-1)})\Phi
((\De\otimes id)(F^{(n-1)})^{-1})((F^{(n-1)})^{-1}\otimes 1)\\
&=1+\xi_n h^n\quad\text{ mod} \quad h^{n+1}.
\endalign$$
 A straightforward calculation shows that if $\Phi$ satisfies
the pentagon identity relative to $\Delta$  then  $\tilde\Phi$
satisfies the pentagon identity, $\widetilde{\Pent}$, relative
to $\tilde\Delta=F^{(n-1)}\Delta (F^{(n-1)})^{-1}.$
The congruence $\Delta=\tilde\Delta$ mod $h$, implies that the
 operator, $\tilde\de$, defined by (30) with  $\tilde\Delta$
replacing $\De$ is congruent mod $h$ to $\de$, therefore
$$1=\widetilde{\Pent}(1 +\xi_n h^n))=1+\tilde\de \xi_n h^n=
1+\de \xi_n h^n\quad\text{mod}\quad h^{n+1}$$
so $\xi_n$ is a $\de$ cocycle.

The second statement follows from the equation
$$B(F^{(n-1)}+ f_n h^n,\Phi)=B(F^{(n-1)},\Phi)+
(\de f_n) h^n\quad\text{mod}\quad h^{n+1}.\qed$$

Since $H^n(\cC)=\wedge^n\cG$, every n-cocycle is
 cohomologous to an element of $\wedge^n\cG$
which we take as the canonical representative of its cohomology
class. Furthermore antisymmetrization annihilates
coboundaries and projects  a cocycle onto
the canonical representative.
Antisymmetrizing  (42)
gives $-\frac23$ times the Yang-Baxter expression $YB(f)$,
$$YB(f):=[f^{12}, f^{13}+ f^{23}] +[f^{13}, f^{23}]=-[f^{23}, f^{12}+
f^{13}] +[f^{12}, f^{13}].\tag44$$

\proclaim{Lemma 3.2}
The necessary and sufficient condition for the existence of a
solution of (37) up to order $h^3$, is that the infinitesimal
$f$ satisfy
$$YB(f)=\frac32\varphi.\tag45$$
\endproclaim
\par\noindent{\bf Proof.}
 Assuming that $f$ satisfies (45), the antisymmetrization of
(42) is $-\varphi$ as required for the obstruction to
be a coboundary. Thus there exists an $f_2$ such that
 $F^{(2)}=1 +fh + f_2 h^2$ satisfies (37) to order $h^3$.\qed

In order to continue the process we need to consider conditions
on $F$ which guarantee that the higher obstructions cocycles cobound.
The strategy suggested by Drinfeld's proof of the existence of $\Phi$
is to impose conditions on $F$ that force the
obstructions to lie in a subcomplex with trivial 3-cohomology.

In the case of a general Lie algebra, there are two obvious symmetries that
we can impose on $F$. The first one,  relating to the antipode, was
introduced in equation (21).
 The second condition on $F$ is
$$ F^{21}_h= F_{-h},\quad\text{or in components}\quad
f_n^{21}=(-1)^nf_n. \tag46$$
This can be interpreted in terms of the $*_h$ product on the function
algebra as saying that the commutator and anticommutator
  are, respectively, expansions in odd and even powers of $h.$

The following lemma shows that in solving equation (37) recursively
using  obstruction  theory we can reduce to a subcomplex of invariants
with respective the symmetries ($\theta,S$) of $\Phi$.

\proclaim{Lemma 3.3 }
Let $\Phi=\Phi_h$ be a solution to the pentagon identity satisfying
(11),(12) and  $\Phi_{-h}=\Phi_h$.
 Let
$$F'=1+f h +\sum_{k=1}^{ n-1}f_k h^k$$
be a solution of  (37) satisfying (46) and (21), all this modulo
$h^n,$ for $n\geq 3.$

(a)  Define $\eta$, $F$ and $\xi$ by
$$
\align(F^{\prime 21})^S F'&=1+ \eta h^n \quad \text{mod}\quad h^{n+1}\\
F&=F' +\frac12\eta h^n\\
B(F,\Phi)&=\xi h^n \quad \text{mod}\quad h^{n+1}.
\endalign
$$
Then $F$ satisfies (21) and (46) mod $h^{n+1}$ and
$\xi^\tau=\xi^S=(-1)^{n+1}\xi,$ where $\tau$ is defined by (31).
 For $n$ odd the obstruction cocycle cobounds and
we can extend $F^{(n-1)}$ to  a solution to (37) to order $h^{n+1}$
which still satisfies  (21) and (46) to the same order.

(b) Let $\theta$ be an involutive automorphism of $\cG$ and
use the same symbol to denote its extension to the complex $\cC$.
Assume that a solution, $\Phi$, to the pentagon equation has been constructed
such that
$\Phi^\theta=\Phi.$
 If $F^\theta=F$ then $\xi^\theta=\xi$. If
$F^\theta_h= F_{-h}$ then $\xi^\theta=(-1)^n \xi.$

\endproclaim

\par\noindent{\bf Proof.}
Define $\hat F_h=F_{-h}.$ Using, successively, the  symmetry (46) on $F$
together with the identities
$\hat\Phi=\Phi,\quad\Phi\Phi^{321}=1$, $\Phi=1$ mod $h^2$ we get,
 modulo $h^{n+1},$
$$
\align
\xi_n^{321}h^n&=B(F, \Phi)^{321}=F^{21}(\De\ot id)(F^{21}) \Phi^{321}
-F^{32}(id\ot \De)(F^{21})\\
&=\hat F^{12}(\De\ot id)\hat F \Phi^{321}-
\hat F^{23}(id\ot \De)(\hat F)\\
&=-B(\hat F,\Phi)\Phi^{321}=-B(\hat F,\Phi)=-\xi_n (-h)^n\\
&=(-1)^{n+1}\xi_n h^n.
\endalign
$$
This shows that the obstruction cocycle at odd order lies in
the acyclic subcomplex $\cC_{\tau, +}$ and so it is a coboundary.
To prove the second identity on $\xi$ in part (a)
we use the additional facts that $F=1$ mod $h$,
 (21) has been solved modulo $h^{n+1}$
 and $(\Phi^{321})^S=\Phi.$ Again, computing modulo
 $h^{n+1},$
$$\align
 (\xi_n^{321})^S h^n&= (B(F, \Phi)^{321})^S=
(\Phi^{321})^S(\De\ot id)(F^{21})^S(F^{21})^{S}
-(id\ot \De)(F^{21})^S(F^{32})^S\\
&=\Phi(\De\ot id)(F^{-1})(F^{12})^{-1}
-(id\ot \De)(F^{-1})(F^{23})^{-1}\\
&=(id\ot \De)(F^{-1})(F^{23})^{-1}B(F,\Phi)(\De\ot
id)(F^{-1})(F^{12})^{-1}\\
&=(id\ot \De)(F^{-1})(F^{23})^{-1}(\xi_nh^n)(\De\ot
id)(F^{-1})(F^{12})^{-1}=\xi_n h^n.
\endalign
$$

Condition (b)  follows from
$$\xi^\theta h^n=B(F,\Phi)^\theta=B(F^\theta,\Phi^\theta)=
B(F,\Phi)=\xi h^n\quad \text{mod}\quad h^{n+1}.$$
The second part is proved in the same way
$$\xi(-h)^n=B(\hat F,\Phi)=B(F^\theta,\Phi^\theta)=B(F,\Phi)^\theta=
\xi^\theta h^n \quad \text{mod}\quad h^{n+1}.$$\qed

\noindent{\bf Remark 3.1} We shall use the concept of a bialgebra action to
 consider invariance with respect to the Cartan
subalgebra, $\cH$,  of $\cG$.  If  $({\Cal U}, \Delta)$ and
$({\Cal P}, \bar\Delta) $ are two  bialgebras, then a bialgebra
action of ${\Cal P}$ on  ${\Cal U}$ is a ${\Cal P}$ module structure on  ${\Cal
U}$  satisfying the two conditions

 $$p\cdot(uv)=\sum (p_{(1)}\cdot u)(p_{(2)}\cdot v),\quad\text{and}
\quad\Delta(p\cdot u)= p\cdot \Delta(u).$$

In item (b) of Lemma 3.1 we have extended
 the action to $\cC^n={\Cal U}^{\ot n}$   by
$$p\cdot (u_1\ot\cdots \ot u_n)=\sum (p_{(1)}\cdot u_1)\ot \cdots \ot
(p_{(n)}\cdot u_n),$$
where $\sum p_{(1)}\ot \cdots \ot p_{(n)}$ represents the $n^{th}$
iteration of $\bar\Delta$.
We say that an element of $\Psi\in\cC$ is invariant under ${\Cal P}$ if
$ p\cdot \Psi=\epsilon(p)\Psi $ for all $p\in {\Cal P}.$
The first of the two properties implies that the product of invariants is an
invariant, thus if $\Phi$ and $F$ are ${\Cal P}$ invariant then so is
$\xi.$
The second property implies that the coboundary
commutes with the action so the invariants form a subcomplex.
 If $\theta$ is the Cartan involution,
 the actions of ${\Cal P}, \theta, S$
are all compatible and we can add the condition of ${\Cal P}$
invariance in Lemma 3.3:

For the Drinfeld-Jimbo infinitesimal, the relevant bialgebra is
${\Cal P}=U(\cH)$  and the relevant involution is
the Cartan involution.  The subcomplex of invariants is a direct summand of
the total complex. The even obstructions are invariant under
$\theta$ and $\cH$. Unfortunately, this  still
does not make the obstructions cohomology classes zero.

The next idea is to cancel the nonvanishing
obstruction cohomology class  using the Chevalley-Eilenberg
cohomology of the Lie
algebra structure on $\cG^*_f$ defined by the infinitesimal $f$.
The fundamental proposition is once again due to Drinfeld.

\proclaim{Proposition 3.1 (Drinfeld\cite{Dr3})}
For $f\in\wedge^2\cG$ the $\cG$ invariance of $YB(f)$ is equivalent to the
Jacobi identity for the bracket $[\lambda,\mu]$ on $\cG^*$
defined by
$$\la [X\ot 1 + 1\ot X, f],\lambda\wedge \mu\ra=\la
X,[\lambda,\mu]\ra.$$
\endproclaim

Now suppose that we are given  $F$ defined  to order $h^{2m-1}$
 with obstruction cocycle $\xi$. If we change $F$
by adding   $\chi\in \wedge^2\cG$
 $$F'=F +\chi h^{2m-1},$$
the new obstruction cocycle is
$$ B(F', \Phi)=(\xi +b(f, \chi))h^{2m},$$
where
$$ b(f,\chi)=f^{23}(\chi^{12}+\chi^{13})
+\chi^{23}(f^{12}+f^{13})-f^{12}(\chi^{13}+\chi^{23})
-\chi^{12}(f^{13}+f^{23}).$$

Projecting onto cohomology by antisymmetrization, we get
$$\Alt(\xi + b(f, \chi))=\Alt (\xi)-
\frac23\widetilde{YB}(f,\chi)\tag47$$
where $\widetilde{YB}$ is the polarization of (42):
$$[f^{12}, \chi^{13}+ \chi^{23}] +[f^{13}, \chi^{23}]+
[\chi^{12}, f^{13}+ f^{23}] +[\chi^{13}, f^{23}].\tag48$$
Thus we need to choose $\chi$ so that (47) is zero.
Formula (48) is a particular example of the {\bf Schouten bracket}
$$\align
[[\cdot,\cdot]]:\wedge^k\cG\ot\wedge^l\cG&
\longrightarrow\wedge^{k+l-1}\cG\\
[[X_1\wedge\cdots\wedge X_k, Y_1\wedge\cdots\wedge Y_l]]&=\sum
(-1)^{i+j}[X_i,Y_j]\wedge X_1\wedge\cdots \hat X_i \cdots
\hat Y_j\cdots \wedge Y_l.\tag49\endalign
$$
Thus (47) is zero if and only if
$$ \frac32 \Alt \xi=[[ f, \chi]].\tag50$$

We remind the reader of the  well known fact.
\proclaim{Proposition 3.2}
If we identify $\wedge^n\cG$ with a skew symmetric multilinear maps on
$\cG^*$, then the Chevalley-Eilenberg
coboundary $d_{CE}$ for the cohomology of the Lie algebra $\cG^*_f$
with coefficients in the ground field is given up to a normalizing factor
by the Schouten bracket, $d_{CE}=[[ f, \cdot].$
\endproclaim
\par\noindent{\bf Proof.}
For $X\in \cG$ and $\lambda,\mu\in \cG^*$, by definition of the
bracket on $\cG^*$ we have
$$\la [[X, f]],\lambda\wedge\mu\ra=\la [X\ot 1 +1\ot
X,f],\lambda\wedge\mu\ra=
\la X,[\lambda,\mu]\ra=\la d_{CE}X, \lambda\wedge \mu\ra$$
Since both $[[\cdot,f]]$ and $d_{CE}$ are derivations of the exterior
algebra, $\wedge \cG$,  and they agree on the generators, they are
equal.\qed

If  $\Alt \xi$ is a coboundary in the Chevalley-Eilenberg
cohomology then we can adjust the extension at order $h^{2m-1}$
by adding $\chi h^{2m-1}$ and cancel the obstruction in $\de$
cohomology at order $h^{2m}$ . Thus we have a kind of secondary
obstruction theory. The next lemma makes this more precise.
\proclaim{Lemma 3.4 }
Let $\xi$ be the obstruction cocycle for
$F^{(n-1)}$ as  defined in (43), then
$d_{CE}\Alt(\xi)=0.$ If $d_{CE}\,\chi=\frac32\Alt(\xi)$ then
the obstruction cocycle for $F'=F+\chi h^{n-1}$ is a $\de$ coboundary.
\endproclaim
\par\noindent{\bf Proof.}
As in the proof of Lemma 3.1 we use the fact that $1+\xi h^n$
satisfies the $\tilde\Delta$ pentagon identity.  Let $\tilde\de$
be given by formula (30) with $\tilde\Delta$ replacing $\Delta$
and $\eta=\Alt \xi.$ Then
$$
\align
\tilde\Delta(u)&=\Delta(u)+ \Delta_1(u) h=\Delta(u) + [f,\Delta(u)] h,
\quad\text{mod}\quad h^2,\\
\text{Since}\quad1&=\widetilde{\Pent}(1 + \xi h^n)= 1+\tilde\de\xi h^n
\quad\text{mod}\quad h^{n+2},\\
\text{we have}\quad0&=\Alt \tilde\de\xi =\Alt(1\ot \xi -(\tilde\De\ot 1\ot
1)\xi
+(1\ot \tilde\De\ot 1)\xi- ( 1\ot 1\ot \tilde\De)\xi +\xi \ot 1\\
&=(-(\tilde\De_1\ot 1\ot 1)\eta
+(1\ot \tilde\De_1\ot 1)\eta-(1\ot 1\ot \tilde\De_1)\eta)h\\
&=(-[f^{12},\De\ot 1\ot 1(\eta)]+[f^{23}, 1\ot \De \ot 1(\eta)]
-[f^{34}, 1\ot 1\ot \De(\eta)])h\\
& =-[[f,\eta]]h=-d_{CE}(\eta)h\quad\text{mod}\quad h^2.
\endalign
$$
The second part of the lemma follows immediately from (50).
In \cite {Dr3}, Drinfeld pointed out the relevance of the Chevalley-Eilenberg
cohomology of $\cG^*_f$ to the study of bialgebra deformations.
This was developed further by LeComte and Roger \cite{LR}.
Our approach is slightly different since we use both the
Chevalley-Eilenberg cohomology and the Cartier cohomology.
\proclaim{Theorem 3.1}
Given $\cG$ be a Lie algebra over a field of characteristic zero and an
  $f\in \wedge^2\cG$ such that $[[f,f]]$ is $\cG$ invariant.
Let  $\cG^*_f$ be the Lie algebra structure
induced by $f$ on $\cG^*$.
Suppose we are in the situation of Remark 3.1 and
a  bialgebra ${\Cal P}$ acts on $U(\cG)$ and that it preserves $\cG$.
Then both ${\Cal P}$ and $\theta$ act
on the Chevalley-Eilenberg cohomology of $\cG^*_f$.
Assume that $f$ is ${\Cal P}$ and $\theta$ invariant.
 If, in  $H^3_{CE}(\cG^*_f)$, the subspace, $H^3_{CE}(\cG^*_f)'$
 of invariants of both actions
is zero,  then there exists a quantized universal
enveloping algebra having undeformed multiplication and deformed
comultiplication with infinitesimal induced by the commutator $[f,\De]$.

\endproclaim
\par\noindent{\bf Proof.}
Define $\varphi=\frac23[[f,f]]$. Then applying Theorem 2.1 we construct
$\Phi$ satisfying (11), (12) and $\Phi^\theta=\Phi$. If there exists a solution
$F_h$ to (26) we can define $\De_h=F_h\De F_h^{-1}.$
Lemma 3.2 says that we can find such a solution modulo $h^3$.
As noted earlier, Lemma 3.3  implies that the only obstructions
lie in the subcomplex of $\theta$ and ${\Cal P}$ invariants and
Lemma 3.4  shows that it is enough to consider the invariants in the
Chevalley-Eilenberg cohomology. If this space of invariants is zero,
there are no obstructions to the recursive solution of (26).\qed

Associated to any deformation of comultiplication in $U(\cG)$
there is a dual deformation of the multiplication in $\cO$, which is expressed
by
$$ f_{m,\lambda}\, *_h \, f_{n,\mu}= f_{m,\lambda}\ot f_{n,\mu}
\circ F_h\Delta F_h^{-1}.\tag51$$
We shall check the compatibility with  the formula
$$ f_{m,\lambda}\, *_h \, f_{n,\mu}= f_{m\tilde\ot n,\mu\tilde\ot
\lambda}.\tag52$$
By the form of the equivalence, equation (20),
$m\tilde\ot n=(F^{-1})_{(1)}(m)\ot (F^{-1})_{(2)}(n)$
and $\mu\tilde\ot\lambda=(F^{-1})_{(1)}(\mu)\ot (F^{-1})_{(2)}
(\lambda)=\mu\circ S((F^{-1})_{(1)})\ot \lambda\circ
S((F^{-1}))_{(2)}=\mu\circ F_{(2)}\ot \lambda\circ F_{(1)},$
where the  last equation follows from (21).  These identities,
 when substituted in the right side of   (51) give (52).
The Hopf algebra $\cO[[h]]$ with this product $*_h$,
undeformed comultiplication, and dual antipode, $\tilde S'$, is the
function algebra of a quantum group  in  a nonstandard
presentation.

\subheading{\S4 Drinfeld-Jimbo quantization of $U(\cG)$}

To construct the Drinfeld-Jimbo QUE algebra
 for a semisimple Lie algebra $\cG$,
choose a Cartan subalgebra and a system of positive roots, $\Pi$.
Let
$$f=\sum_{\alpha\in \Pi}X_{\alpha}\wedge X_{-\alpha}\tag51$$
 be the Drinfeld-Jimbo classical $R$-matrix. It is
$\cH$ invariant and skew invariant under the Cartan involution.
We shall apply Theorem 3.1 to the case when ${\Cal P}=U(\cH)$ and
$\theta$ is the Cartan involution.

If we use the Killing form on $\cG$ to define an isomorphism between
$\cG$ to $\cG^*$ then the Lie algebra structure $\cG^*_f$ can be
transferred  to a new Lie algebra structure, denoted $\cG_f$,
 on $\cG.$
\proclaim{Lemma 4.1}
The space of invariants under $\cH$ and $\theta$, $H^3_{CE}(\cG_f)'$,
is trivial.
\endproclaim
\par\noindent{\bf Proof.}
The  Lie algebra, $\cG_f$, can be described quite simply
in two steps. Let $\cN_\pm$ be the two subalgebras of $\cG$ consisting of
the sum of the positive root spaces and
the sum of the negative root spaces, respectively.
 First form the Lie algebra sum $\cN_+\oplus \cN_-$
then take the semidirect product with $\cH$ where $\cH$ acts on
$\cN_+$ by the standard bracket and on $\cN_-$  by the negative of
the standard bracket.

In general, given an abelian subalgebra
acting semisimply under the adjoint representation, the
non-trivial cohomology lies entirely in the subcomplex of zero
weight. For the Lie algebra $\cG_f$  the relevant subcomplex
of $\wedge \cG$ is just $\wedge \cH$. Since $f$ is $\cH$ invariant
the Schouten bracket $[[f,\cdot]]$ restricts to the zero operator
on the subcomplex
 and therefore the cohomology is
$H^n_{CE}(\cG_f^*)=\wedge^n\cH.$
The Cartan involution restricted to $\cH$ is $-1$ so the
subspace of invariants in any odd exterior power is zero.
\qed

We can now state our main result which follows immediately from
Theorem 3.1.

\proclaim{Theorem 4.1}
Let $\cG$ be a semisimple Lie algebra over a field of characteristic
zero and $f$ the Drinfeld-Jimbo classical $R$-matrix as defined in
(51). Then there exists an $F=F_h\in U(\cG)^{\ot 2}[[h]]$ such that
\roster
\item
$F$ transforms the associativity constraint to the identity (equation (26)).
\item
$F$ is invariant under $\cH$ and any automorphism group of $\cG$ under
which $f$ is invariant.
\item
$F^{\theta}_h=F^{21}_h=F_{-h},$ where $\theta$ is the Cartan
involution.
\item
$F^SF^{21}=1,$ where $S$ is the antipode.
\endroster
\endproclaim

\proclaim{Corollary 4.1}
Let $\cG$ and  $F$ be as in Theorem 4.1.
Define a quasi-bialgebra deformation of $U(\cG)$ by leaving the
multiplication undeformed and deforming the comultiplication by
twisting by $F_h$
$$\Delta_h(u)=F_h\Delta(u) F_h^{-1}.$$
The resulting deformation has the following properties:
\roster
\item
The deformed comultiplication is coassociative.
\item
The restriction of $\Delta_h$ to $U(\cH)$ is the standard, undeformed,  so
$U(\cH)$ is undeformed.
comultiplication.
\item
$\theta$ is a coalgebra antiautomorphism relative to $\Delta_h$.
\item
The undeformed antipode $S$ is a coalgebra antiautomorphism relative
to $\Delta_h$.
\item
Let $F=\sum F_{1i}\ot F_{2i}$ and $w=\sum F_{2i}S(F_{1i}).$ Then
$$\tilde S(u)=w^{-1}S(u) w$$ defines an antipode relative to which
the deformation is a Hopf algebra.
\endroster
\endproclaim
\par\noindent{\bf Proof.}
Each item, except the last, follows from the corresponding item in Theorem 4.1.
The last statement follows from the explanation of the twisting
transformation and equation (27) with $\alpha=1$.\qed

Regarding the uniqueness of the solution, $F$, with a given
infinitesimal, $f$, we have the following proposition.
\proclaim{Proposition 4.1}
Let $\cG$ be  a Lie algebra of a field of characterisitic zero,
and $\Phi\in U(\cG)^{\ot 3}[[h]]$ a solution to the pentagon identity.
Given two solutions, $F_h$ and $F'_h$, to equations (21) and (26)
both with initial term $1\ot 1$ and
the same infinitesimal $f\in\wedge^2\cG$, then there
exists a $u_h\in U(\cG)[[h]]$ such that
$$ F_h=(u_h \ot u_h) F'_h \De(u_h^{-1}),\quad\text{and}\quad u_hS(u_h)=1.$$
\endproclaim
\par\noindent{\bf Proof.}
 We shall  prove, as usual, that if $F$ and $F'$ agree to order $h^n$
then there exists a $u=1 + u_n h^n$ such that mod $h^{n+1}$ $F$ equals the
transform of $F'$,  $(u\ot u )F'
\De(u^{-1})=F'+\de(u_n)h^n$.
The facts that $F$ and $F'$ both satisfy (26) and that they
agree to order $h^n$ imply that the difference of the
 coefficients of $h^n$  is a $\de$ cocycle, $\de(f_n-f'_n)=0$.
Therefore, there exists a $v$ such that $f_n-f'_n=\de v + w$,
where $w\in \wedge^2 \cG$. However (21) implies that $f_n + f_n^S=
f'_n + (f'_n)^S.$ Thus $(f_n-f'_n)=-(f_n -f'_n)^S.$ Now $w=S(w)$ and
$(\de v)^S=\de(S(v))$, so we have
$f_n-f'_n=\frac12((f_n-f'_n)-(f_n -f'_n)^S)=\de(\frac12(v -S(v)).$
Set $u_n=v-S(v)$.

For an arbitrary pair, $F_h$ and $F'_h$, the  $h$ expansion of
$u_h$ is given by an infinite product\newline
 $Lim_{n\rightarrow\infty}(1+u_nh^n)\cdots(1+u_2h^2)$.  Next
consider the condition
$S(u_h)u_h=1.$ Suppose that we have defined $u'$
which is the product of the first $n$ terms, so that, modulo $h^{n+1}$,
 $F$ and $F'$ agree and that $S(u')u'=1 +\xi_n h^n$.
 The fact that both $F=(u'\ot u')F'(\De(u')^{-1})$ and $F'$ satisfy
(21) modulo $h^{n+1}$ implies that  ( modulo $h^{n+1}$ )
$$\De(S(u')^{-1})(F')^S (S(u')\ot S(u'))(u'\ot u')F'(\De(u')^{-1})=1.$$
Substituting in $(1 +\xi_n h^n)u^{\prime-1}$ for $S(u')$
and expanding the product, using (21) and $F=F'=1$ modulo $h$,
we see that $\xi_n\in U(\cG)$ is  a $\de$ cocycle, i.e., primitive
relative to $\De$, so it is an element of $\cG$.  From the
definition it also follows that $\xi$ is $S$ invariant, so, in fact, it
is zero and no correction is necessary.\qed

In Theorem 4.1 we allow $f$ to be any linear combination of the
Drinfeld-Jimbo infinitesimals for the simple factors.
This possibility is important for applications to homogeneous spaces,
in particular the symmetric space $(G\times G)/G.$
When $\cG$ is simple, the  deformation described in Theorem 4.1
is unique up to inner automorphism and
change of parameter, as was noted in \cite{Dr3}. For a complete proof
see \cite{SS}.

In this way it is possible to quantize some of the classical
$R$-matrices classified by  Belavin Drinfeld. These examples will
be discussed in a future paper.

Another interesting application of Theorem 3.1
is the quantization of  an arbitrary infinitesimal $f$ with
invariant Schouten bracket $[[f,f]]$ for any three dimensional
Lie algebra. In this case $\wedge^3\cG$ is one dimensional, if
$[[f,\chi]]\neq 0$ for some $\chi\in\wedge^2\cG$
then $H^3_{CE}(\cG^*_f)=0.$ On the other hand, suppose $[[f,\chi]]=0$ for
all $\chi$. We can choose a basis $\{ X, Y, Z\}$
such that  $f=X\wedge Y$. The condition that $[[f,f]]=0$ implies
that $[X,Y]\in$ span$\{X,Y\}$. The conditions  $[[f,Y\wedge Z]]=0$
and $[[f, X\wedge Z]]=0$
imply that $[X,Y]\in$ span $\{Y,Z\}$ and
 $[X,Y]\in$ span $\{ X, Z\}$ respectively. All the conditions together
imply $[X,Y]=0$. In this case $F=e^{hf}$ satisfies
(26) with $\Phi=1$ and gives the desired quantization.

\subheading{\S Appendix. Auxiliary conditions on the $R$-matrix}

In this appendix we return to the braiding and consider the problem of proving
the
existence of a quasi-triangular quasi-Hopf deformation
 $(U(\cG)[[h]],\cR_h,\Phi_h)$
which includes auxiliary conditions on $\cR$  corresponding to the
conditions on $\Phi$ appearing in Theorem 2.1.
The correspondence is shown in the table given in the following
proposition, where
$S$ is the antipode, $(\hat\cR)_h=\cR_{-h}$ and similarly for $\Phi.$
\proclaim{Proposition}
For any $\cG$ invariant symmetric invariant element  $t\in\cG\ot\cG$ and
automorphism $\theta$
preserving $t$
there exists a pair of $\cG$ invariant elements
 $\Phi\in U(\cG)^{\ot 3}[[h]]$ and $\cR\in
 U(\cG)^{\ot 2}[[h]]$,  where $\cR\equiv 1 +h t$ mod $h^2$,
which satisfy the pentagon and hexagon identities
and have  the following symmetries.
$$
\align
\Phi^{321}\Phi=1 &\hskip.5cm  \cR^{21}=\cR,\tag A.1\\
\Phi^{\theta}=\Phi &\hskip.5cm \cR^{\theta}=\cR,\tag A.2 \\
\Phi^S\Phi=1 &\hskip.5cm \cR^{S}=\cR,\tag A.3\\
\hat{\Phi}=\Phi&\hskip.5cm  \hat {\cR} \cR=1.\tag A.4
\endalign
$$
\endproclaim
\noindent {\bf Proof}
 In \cite {Dr2} Drinfeld proved the existence of a pair $(\Phi,\cR)$
satisfying the pentagon and hexagon identities using standard Cartier
coalgebra coboundary, $\delta$, to study the pentagon equation and a modified
Cartier coboundary, $\delta'$, in which the last factor is  ``frozen'' to study
the hexagon identities.  In the proof he  imposes the condition (A.1).
In fact, the remaining conditions can be included in his proof. We
shall prove this by showing that the crucial obstruction
equations can be solved in the appropriate subcomplex.

As usual, suppose that we have a pair $(\Phi,\cR)$
giving a truncated deformation  defined to order $h^n$
and we want to extend it to order $h^{n+1}.$
(All further equations which contain $h^n$ will be understood to be modulo
$h^{n+1}.$)
Define the
obstruction pair $(\xi, \psi)$ by
$$\align \Pent(\Phi)&=1 +\xi h^n \quad \text{mod}\quad h^{n+1}\tag A.5\\
(\De\ot id)\cR&=\Phi^{312}\cR^{13}(\Phi^{132})^{-1}\cR^{23}\Phi
+\psi h^n  \quad \text{mod}\quad h^{n+1} \tag A.6
\endalign
$$

Assuming $\cR$ symmetric and $\Phi^{321}\Phi=1$ mod $h^{n+1}$,
transposing tensor factors 1 and 3, (A.6) gives
the other hexagon identity with error term
(obstruction cochain) $\psi^{321}$.

By the discussion
in \S 2 we know that there is  an extension
$$\Phi'=\Phi + \phi h^n,$$
such that the pentagon identity together with  all the
identities (A.1-3) are satisfied to order $h^{n+1}$.
The element $\phi$ can be modified by a $\delta$ cocycle with
 the appropriate symmetries. The condition  (A.4), which did not
appear in \S2, is trivial at this stage, since we can simply assume
that the deformation parameter is $h^2$. This means that the
obstruction $\xi$ is automatically zero at  odd orders. However,
when we turn to the hexagon identities, the obstruction $\psi$
is not necessarily zero at odd order and Drinfeld uses  the freedom of
adding a $\delta$ cocycle to $\Phi$ in ``killing the obstruction''.
We must show that when we impose the extra symmetries, no modification
of $\Phi$ is necessary at odd orders. More explicitly, if we extend
$\Phi$ as above and $\cR$ by
$$\cR'=\cR + r h^n,$$
then the new  obstruction is
$$\psi -(\De \ot id)r +\phi^{312} + r^{13} -\phi^{132} + r^{23}
+\phi.$$
The terms in $r$ define the modified Cartier operator $\delta'r
=(\De \ot id) r - r^{13} -r^{23}$.
Transposing factors 1,2 and subtracting gives
$$\psi -\psi^{213} + 6 \operatorname{Alt}\phi.$$
Drinfeld proves that $\psi-\psi^{213}\in\wedge^3\cG$, and  it is
possible to cancel this term by adding a $\delta$ cocycle to $\Phi$
(since $\delta(\wedge^3\cG)=0$).
We want to show that for $n$ odd, $\psi-\psi^{213}=0.$
At this point the additional conditions (A.3-4) on $\cR$ become
relevant. Suppose that these conditions are satisfied modulo $h^n$,
and define $\chi$ by
$$ \hat\cR\cR=1 + \chi h^n.\tag A.7$$
Since $\cR=1$ mod $h$ we also have
$$\cR\hat\cR=1+ \chi h^n  .\tag A.8$$
Substituting $-h$ for $h$ we find
$$\hat \cR\cR= 1+ \chi (-h)^n .$$
Therefore, $\chi=0$ for $n$ odd. For $n$ even, applying  $S$ to (A.8), we get
$$(\cR)^S(\hat\cR)^S=1+\chi^S h^n,$$
so $\chi^S=\chi.$ Setting
$$\cR'=\cR - \frac12 \chi h^n,$$
we have (A.3-4) to order $h^{n+1}.$ We compute the obstruction $\psi$
for this  $\cR'$, which we shall denote simply $\cR$.
Applying $S$,
$$\align
(\De\ot id)\cR=(\De\ot id)(\cR)^S&=(\Phi^{312}\cR^{13}
(\Phi^{132})^{-1}\cR^{23}\Phi)^S+\psi^S h^n \\
&=\Phi^{S}(\cR^{23})^{S}((\Phi^{132})^{-1})^{S}(\cR^{13})^{S}
((\Phi^{312})^{S}
+\psi^S h^n\\
&=\Phi^{-1}\cR^{23}\Phi^{132}\cR^{13}(\Phi^{312})^{-1}
+\psi^S h^n\\
&=\Phi^{321}\cR^{23}(\Phi^{231})^{-1}\cR^{13}\Phi^{213}
+\psi^S h^n,  \endalign
$$
where we used $\Phi^{321}=\Phi^{-1}$ in the last step.
Transposing factors 1,2 leaves the left side invariant and gives
$$ (\De \ot id)R =\Phi^{312}\cR^{13}
(\Phi^{132})^{-1}\cR^{23}\Phi + (\psi^{213})^S h^n.$$
Therefore
$$(\psi^{213})^S =\psi.\tag A.9$$

On the other hand
$$\align
((\De\ot id)\hat{\cR})^S&=(\hat\Phi^{312}\hat{\cR}^{13}
(\hat\Phi^{132})^{-1}\hat{\cR}^{23}\hat\Phi)^S
+\psi^S (-h)^n \\
&=\Phi^{-1}(\cR^{23})^{-1}\Phi^{132}(\cR^{13})^{-1}(\Phi^{312})
+(-1)^n\psi^S h^n.
\endalign$$

Taking inverses in (A.6) gives
$$
\align
((\De\ot id)\hat\cR)^S&=(\De\ot id)\cR^{-1}=(\Phi^{312}\cR^{13}
(\Phi^{132})^{-1}\cR^{23}\Phi+\psi h^n)^{-1} \\
&=\Phi^{-1}(\cR^{23})^{-1}\Phi^{132}(\cR^{13})^{-1}(\Phi^{312})
-\psi h^n.
\endalign
$$
Thus
$$\psi^S=(-1)^{n+1}\psi.\tag A.10 $$
Then
$$\psi^{213}=\psi^S=(-1)^{n+1}\psi,$$
so $\psi$ is symmetric in 1,2 for $n$ odd, and antisymmetric in 1,2 for
$n$ even.

For $n$ even we can ``kill the obstruction'' by adding $\phi h^n$ to $\Phi$
where $\phi=\frac13 \psi\in \wedge^3 \cG$.
For $n$ odd, we don't change  $\Phi$ but add to $\cR$ a term $r
h^n$ satisfying the equation
$$\delta' r=\psi,\tag A.11$$
where $\delta'$ is the modified Cartier coboundary in dimension 2.
The second cohomology group for $\delta'$ is $\wedge^2\cG\ot U(\cG)$,
and since $\psi$ is symmetric in 1,2 it must cobound, moreover
$\psi$ is $S$ invariant and the cochain it cobounds can be chosen $S$
invariant. Thus the condition (A.3) is preserved.  Moreover, since
$n$ is odd, and $\cR=1$ mod $h$,  adding any term $r h^n$ will not
affect condition (A.4). This shows that one can take the next step in
the recursive construction of the pair $(\Phi, \cR)$.\qed

Consider the regular representation of $U(\cG)$ on the
 the compactly supported functions
on the corresponding Lie group $G$. If we define a scalar product using
the Haar measure, then for $u\in U(\cG)$
 the operator $S(u)$ is a  formal
adjoint to the operator $u$. The elements $F, R$ of the
double tensor product and $\Phi$ of the triple tensor product
can be considered as formal power series operators acting on the
functions on $G\times G$ and $G\times G\times G$ respectively.
Extending $S$ conjugate linearly and using the pure imaginary
deformation parameter $i\nu$ conditions (A.3) and (A.4)  and the
conditions (21) and (46)  on $F$ imply  formal unitarity.
\proclaim{Corollary}
Extending $S$  conjugate linearly and using it to define the formal
adjoint, as explained above, we have the following formal unitarity
conditions on $F, \cR$ and $\Phi$:
$$\align
F_{i\nu}F_{i\nu}^*=&F_{i\nu}F_{-i\nu}^S=1\\
(\Phi_{i\nu})^*\Phi_{i\nu}=&\Phi^S_{-i\nu}\Phi_{i\nu}=1
\\
(\cR_{i\nu})^*\cR_{i\nu}=&\cR^S_{-i\nu}\cR_{i\nu}=1.
\endalign
$$
\endproclaim

\end